\documentclass[11pt]{article}
\usepackage[a4paper,margin=1in]{geometry}
\usepackage{lmodern}
\usepackage[round,sort&compress]{natbib}

\usepackage[utf8]{inputenc}
\usepackage[T1]{fontenc}
\usepackage[english]{babel}

\usepackage{graphicx}
\usepackage{booktabs}
\usepackage{makecell}
\usepackage{array}
\usepackage{multirow}
\usepackage{float}
\usepackage{subcaption}
\usepackage{wrapfig}

\usepackage{amsmath,amssymb,amsfonts,amsthm}
\usepackage{mathtools}
\usepackage{dsfont}
\usepackage{bm}
\usepackage{nicefrac}
\usepackage{microtype}

\usepackage{algorithm}
\usepackage{algpseudocode}

\usepackage{enumitem}
\setlist{noitemsep,topsep=2pt,parsep=0pt,partopsep=0pt}
\usepackage{xcolor}
\usepackage{url}
\usepackage[hidelinks]{hyperref}
\usepackage[nameinlink,capitalise]{cleveref}
\usepackage{verbatim}   % \begin{comment} ... \end{comment}

\theoremstyle{plain}
\newtheorem{theorem}{Theorem}
\newtheorem{proposition}{Proposition}

\theoremstyle{definition}
\newtheorem{problem}{Problem}
\newtheorem{definition}{Definition}

\theoremstyle{remark}
\newtheorem{remark}{Remark}

\crefname{problem}{Problem}{Problems}
\Crefname{problem}{Problem}{Problems}

\newcommand{\R}{\mathbb{R}}
\newcommand{\N}{\mathbb{N}}

\newcommand{\de}{\mathrm{d}}
\newcommand{\norm}[1]{\left\lVert#1\right\rVert}
\DeclareMathOperator{\supp}{supp}
\DeclareMathOperator{\dist}{dist}
\DeclareMathOperator{\disp}{disp}

\DeclareMathOperator{\vol}{vol}

\newcommand{\dH}{d_{\mathrm{H}}}

\title{Recovering Lower-Dimensional Semialgebraic Support of a Measure from its Moments}
\author{Ruben Karapetyan, Shenyuan Ma, Ale\v{s} Wodecki, Jakub Mare\v{c}ek%
\thanks{Corresponding author: \texttt{jakub.marecek@fel.cvut.cz}}\\[6pt]
Czech Technical University in Prague}
\date{}
\begin{document}

\maketitle

\begin{abstract}
Recovering probability measures from their moments has numerous applications, esp. in connection with the method of moments in statistics and optimization. In the setting where measure need not be finitely atomic, but its support is known to be compact and semialgebraic with codimension at least one, the problem is still open. We combine moment-matrix kernel information with the Christoffel--Darboux kernel to provide a discrete approximation of the support. To validate the proposed approach, we test our algorithm on analytically computed moments and pseudo-moments arising from polynomial optimization problems without unique global minimizers. This complements well-known recent work on recovery of measures with algebraic support, where the kernel of a moment matrix can reveal polynomials vanishing on the support, and on recovery of sufficiently regular full-dimensional supports, where estimators constructed by thresholding the Christoffel--Darboux kernel are known to converge asymptotically to the support.
\end{abstract}

\section{Introduction}
\label{sec:intro}

The recovery of a measure from finitely many of its moments is a classical
problem \citep{akhiezer1965classical,Schmudgen2017}. It underlies the method of
moments in statistics \citep{lindsay1995} and the extraction of solutions in
polynomial optimization \citep{Lasserre2001GlobalOptimization,
Parrilo2003SemidefiniteProgramming}. In the second of these settings one does
not, strictly speaking, observe moments at all. The moment relaxation of degree
$d$ of a polynomial optimization problem (POP) returns a \emph{pseudo-moment}
sequence, that is, the values of a linear functional which behaves on
polynomials of low degree like integration against a probability measure, but
which need not be the moment sequence of any measure
\citep{lasserre2013moment,laurent2009sums}. The optimal values of the
relaxations converge to that of the POP as the degree increases
\citep{Lasserre2001GlobalOptimization}, but the convergence is asymptotic, and
in practice only small degrees are affordable. One is therefore left with the
problem of recovering the set of minimizers from a pseudo-moment sequence which
is both short and inexact.

Two families of methods are available, and they rest on hypotheses of different
kinds. The first is algebraic and uses the kernel of the moment matrix. Every
polynomial in $\ker M_d(y)$ vanishes on $\supp\mu$, and so the kernel determines
the Zariski closure of the support
\citep{laurent2009sums,wageringel2022truncated,lasserre2013moment}. If $\mu$ is
finitely atomic and the moment matrix admits a flat extension, the atoms can be
computed exactly by solving an eigenvalue problem \citep{CurtoFialkow1998,
HenrionLasserre2005DetectingGlobalOptimality,LopezQuijorna2021DetectingGNS}.
The kernel carries no information about inequalities, however: the semicircle
$\{x_1^2+x_2^2=1,\;x_2\ge0\}$ and the full circle have the same vanishing ideal,
and hence the same moment-matrix kernel in every degree. The second family is
analytic and thresholds the Christoffel--Darboux (CD) kernel. The Christoffel
function is relatively large on the support of the measure and decreases
rapidly away from it, and the sublevel sets of the CD kernel converge in the
Hausdorff metric to the support, provided that the measure is absolutely
continuous and supported on a \emph{full-dimensional} set of sufficient
regularity \citep{lasserre2019empirical,pauwels2021data,vu2022rate,
Lasserre_Pauwels_Putinar_2022}. Inequalities cause no difficulty here, but the
convergence results do not apply once the support has positive codimension,
and this is the usual situation in polynomial optimization, where sets of
minimizers generally have empty interior.

This paper is concerned with the intermediate case: supports which are compact
and semialgebraic, have codimension at least one, and are not necessarily
finite. We are given the moments of $\mu$ up to degree $2d$, but no samples, and
we wish to compute a finite set of points at small Hausdorff distance from
$\supp\mu$. We are not aware of any existing method which does this from
moments alone.

We propose a procedure, \Cref{alg:support}, which uses the moment matrix twice,
first algebraically and then analytically, and we study its behaviour
numerically. Our contributions are as follows.
\begin{enumerate}
\item A method which combines the two techniques (\Cref{sec:method}). The
  algebraic carrier $\mathcal Z=V_\R(\ker M_d(y))$ is recovered first; $\mathcal Z$ is
  then sampled uniformly with respect to its intrinsic volume, by the
  random-slicing method of \citet{Breiding2020}; and only then is the
  \emph{relative} CD kernel thresholded \emph{along} $\mathcal Z$.
  The equations are treated using
  the kernel of the moment matrix and the inequalities using the CD kernel, each
  instrument being used where it is appropriate.
 Since $\kappa_d$ is a polynomial and $\mathcal Z$ is algebraic,
  the sublevel set $\{x\in \mathcal Z:\kappa_d(x)\le\tau_d\}$ is already a semialgebraic
  description of the estimated support, of which the returned cloud of points
  is a sample.
\item Experiments upon exact moments and upon pseudo-moments alike
  (\Cref{sec:analytic,sec:pop}). We work with moments computed in closed form
  for measures satisfying our hypotheses exactly, and with pseudo-moments
  produced by TSSOS \citep{magron2021tssos,wang2021tssos} for POPs whose
  minimizer sets are a segment and a semicircle, where the flat extension
  property fails and the GNS construction is of no help.
\item The ratio of the
 sampling errors
 on the
   ambient grid and the recovered carrier
(with equal budget) grows exponentially with $n$, whatever function
 is filtered on the grid.
  This is shown by  \Cref{thm:grid} and corroborated by experiments.
\end{enumerate}

\section{Related work}
\label{sec:related}

\paragraph{Vanishing ideals from moment matrices.}
That $\ker M_d(y)$ consists of polynomials vanishing on $\supp\mu$ is classical
\citep{laurent2009sums}. To pass from a numerical basis of the kernel to a set
of \emph{generators}, one must separate the new relations in each degree from
polynomial multiples of relations found in lower degrees; border bases,
Gr\"obner-type reductions and real radical computations may all be used for
this purpose \citep{lasserre2013moment,laurent2009generalized}. On truncated
moment problems on varieties and semialgebraic sets see
\citet{wageringel2022truncated,Nie2014_ATKMP}.

\paragraph{Finitely atomic measures.}
When $\mu$ is a finite sum of Dirac measures, the flat extension theorem
\citep{CurtoFialkow1998} gives exact recovery, as implemented in GloptiPoly
\citep{HenrionLasserre2005DetectingGlobalOptimality} and in the truncated GNS
construction \citep{LopezQuijorna2021DetectingGNS}; the extraction is stable
under perturbation \citep{KlepPovhVolcic2018_MinimizerExtractionIsRobust} and
rests on Stickelberger's theorem \citep{cox2020stickelberger}, and the same
ideas appear in super-resolution \citep{Vetterli2002FRI,DeCastro2012BLASSO,
Duval2015SparseSpikes}. None of these methods applies when $\mu$ has a
non-atomic part.

\paragraph{Christoffel--Darboux kernels and the estimation of supports.}
The empirical Christoffel function of \citet{lasserre2019empirical} estimates
the support of a measure from samples; its sublevel sets converge in the
Hausdorff metric under assumptions of regularity and density, and
\citet{vu2022rate} give rates for finite samples. \citet{pauwels2021data}
develop the corresponding theory for moments, including the quotient
construction by which the kernel is defined on a variety, and
\citet{Lasserre_Pauwels_Putinar_2022} give a book-length account; see also
\citet{totik2010christoffel} on asymptotics, \citet{marx2021semialgebraic} on
semialgebraic approximation, and \citet{devonport2021datadriven2} for an
application. All the convergence results known to us assume that the support
is regular closed in the ambient space, $S=\overline{\operatorname{int}S}$, and
in our setting it is not. The mollified CD kernels of
\citet{bentancur2026mollified} are designed for the lower-dimensional case;
we discuss them in \Cref{sec:cd}. For sampling on varieties we use
\citet{Breiding2020} (\Cref{sec:sampling}); on sampling density and topology
see \citet{NiyogiSmaleWeinberger2008}.

\section{Problem setting}
\label{sec:setting}

\subsection{Notation}
\label{sec:notation}

We work throughout in $\R^n$. For $\alpha\in\N^n$ we write
$x^\alpha=x_1^{\alpha_1}\cdots x_n^{\alpha_n}$ and $|\alpha|=\sum_i\alpha_i$.
We denote by $\R[x]_d$ the space of polynomials of degree at most $d$, by
$v_d(x)$ the vector of all monomials of degree at most $d$, and by
$\sigma(d)=\binom{n+d}{d}=\dim\R[x]_d$ the number of such monomials. We
identify $p\in\R[x]_d$ with the vector $\vec p\in\R^{\sigma(d)}$ of its
coefficients in the monomial basis. Given a truncated sequence
$y=(y_\alpha)_{|\alpha|\le 2d}$, the \emph{moment matrix}
$M_d(y)\in\R^{\sigma(d)\times\sigma(d)}$ is defined by
\begin{equation}
\label{eq:moment-matrix}
  [M_d(y)]_{\alpha,\beta}=y_{\alpha+\beta},
  \qquad |\alpha|,|\beta|\le d ,
\end{equation}
so that $\vec p^\top M_d(y)\, \vec q=\int pq\,\de\mu$ whenever $y$ is the
moment sequence of $\mu$. For a set $I$ of polynomials we write
$V_\R(I)=\{x\in\R^n: f(x)=0\ \forall f\in I\}$ for its real zero set, and we
write $\dH(A,B)=\max\{\sup_{a\in A}\inf_{b\in B}\norm{a-b},\ \sup_{b\in
B}\inf_{a\in A}\norm{a-b}\}$ for the Hausdorff distance.

\subsection{Problem statement}
\label{sec:problem}

\begin{problem}
\label{problem_main}
Let $\mu$ be a Borel measure on $\R^n$ supported on a compact semialgebraic set
$S$ given by
\begin{align}
\label{eq:generators}
    & f_i(x) = 0, && i =1,2,\ldots,k, \\
\label{eq:boundaries_given}
    & -1 \leq x_j \leq 1, && j = 1,2,\ldots,n, \\
\label{eq:boundaries}
    & g_j(x)\geq 0, && j =1,2,\ldots,l ,
\end{align}
where \eqref{eq:boundaries_given} normalizes the ambient box. Let
\begin{equation*}
    \mathcal Z:=\left\{x\in\R^n: f_i(x)=0,\ i=1,\ldots,k\right\},
\end{equation*}
suppose that $\mathcal Z$ has dimension $s<n$, and let
$\lambda:=\mathcal H^s|_{\mathcal Z}$ be its intrinsic reference measure, where
$\mathcal H^s$ denotes $s$-dimensional Hausdorff measure. Suppose that $\mu$ is
absolutely continuous with respect to $\lambda$ on $S$,
\begin{equation*}
    \de\mu(x)=w(x)\,\de\lambda(x),\qquad x\in S,
\end{equation*}
with density bounded away from zero, $w(x)\ge w_0>0$ on $S$. Suppose finally
that the moments of $\mu$ up to some degree $2d$ are known,
\begin{equation*}
    y_\alpha=\int x^{\alpha}\,\de\mu(x),
    \qquad \alpha\in\N^n,\ |\alpha|\leq 2d,
\end{equation*}
and that \emph{no} samples from $\mu$ are available. Given $\delta>0$, find a
finite set $S_\delta$ such that
\begin{equation}
\label{eq:goal}
    \dH(S_\delta,\supp \mu) < \delta .
\end{equation}
A set $S_\delta$ satisfying \eqref{eq:goal} is called a
\emph{$\delta$-accurate approximation} of $\supp\mu$.
\end{problem}

Three features of \Cref{problem_main} determine our approach: the support has
empty interior in $\R^n$, so the full-dimensional theory of the CD kernel does
not apply; $S$ is described by equations and inequalities together; and only
moments are given, so any samples must be generated by the method itself
(\Cref{sec:sampling}).

\section{Method}
\label{sec:method}

The procedure is set out in \Cref{alg:support}. The moment matrix is used
twice: its kernel gives the algebraic carrier $\mathcal Z$
(\Cref{sec:kernel_method}), and its pseudoinverse gives the relative CD kernel,
which serves as a membership oracle on $\mathcal Z$ (\Cref{sec:cd}). Between
these two steps, \Cref{sec:sampling} generates the points of $\mathcal Z$ at
which the oracle is evaluated. \Cref{sec:output} describes the output and the
choice of parameters.

\begin{algorithm}[t]
\caption{$\texttt{AlgSup}$: recovery of a lower-dimensional semialgebraic support}
\label{alg:support}
\begin{algorithmic}[1]
\Require Degree $d$; moments $y=(y_\alpha)_{|\alpha|\le 2d}$; sampling budget
         $N$; threshold $\tau_d$.
\Ensure Semialgebraic set $S_{d,\tau_d}$ and a finite sample
        $\mathcal M_{\mathrm{adm}}\subset S_{d,\tau_d}$.
\Statex
\State \textbf{Recover the algebraic carrier $\mathcal Z$.}
  \Statex\hspace{\algorithmicindent}\textbf{(a)} Form the moment matrix $M_d(y)$ as in \eqref{eq:moment-matrix}.
  \Statex\hspace{\algorithmicindent}\textbf{(b)} Compute $I_d:=\ker M_d(y)$ and extract generators degree by degree.
  \Statex\hspace{\algorithmicindent}\textbf{(c)} Set $\mathcal Z:=V_\R(I_d)$.
\Statex
\State \textbf{Sample the carrier.}
  \Statex\hspace{\algorithmicindent}\textbf{(a)} Draw $N$ points on $\mathcal Z$, uniformly with respect to $\mathcal H^s|_{\mathcal Z}$, by random slicing \citep{Breiding2020}.
  \Statex\hspace{\algorithmicindent}\textbf{(b)} Denote the resulting set of points by $\mathcal M$.
\Statex
\State \textbf{Form the relative CD kernel on $\mathcal Z$.}
  \Statex\hspace{\algorithmicindent}\textbf{(a)} Compute the pseudoinverse $A:=M_d(y)^\dagger$.
  \Statex\hspace{\algorithmicindent}\textbf{(b)} Set $\kappa_d(x):=v_d(x)^\top A\,v_d(x)$ for $x\in \mathcal Z$.
\Statex
\State \textbf{Threshold.}
  \Statex\hspace{\algorithmicindent}\textbf{(a)} Evaluate $\kappa_d$ on $\mathcal M$.
  \Statex\hspace{\algorithmicindent}\textbf{(b)} Retain $\mathcal M_{\mathrm{adm}}:=\{x\in\mathcal M:\kappa_d(x)\le\tau_d\}$.
\Statex
\State \Return
$S_{d,\tau_d}=\{x\in \mathcal Z:\kappa_d(x)\le\tau_d\}$ and
$\mathcal M_{\mathrm{adm}}$.
\end{algorithmic}
\end{algorithm}

\subsection{Step 1: the algebraic carrier from the kernel of the moment matrix}
\label{sec:kernel_method}

Let $f_1,\dots,f_k$ be the polynomials in \eqref{eq:generators}, and let
$\mathcal Z$ be their real zero set. The first step computes polynomials in the
vanishing ideal of $\supp\mu$. These need not be the $f_i$ themselves, since
the ideal generated by $f_1,\dots,f_k$ need not be real radical
\citep{lasserre2013moment}. The basis of this step is the following classical
correspondence \citep{laurent2009sums}.

\begin{proposition}[Kernel of the moment matrix and support]
\label{prop:moment-kernel-support}
Let $\mu$ be a finite nonnegative Borel measure on $\R^n$ with moments
$y_\alpha=\int x^\alpha\,\de\mu(x)$ up to degree $2d$, and let $M_d(y)$ be the
associated moment matrix. Then for every $f\in\R[x]_d$,
\begin{equation*}
    f\in\ker M_d(y)
    \quad\Longleftrightarrow\quad
    \supp\mu\subseteq V_\R(f):=\{x\in\R^n: f(x)=0\} ,
\end{equation*}
or equivalently
$\ker M_d(y)=\{f\in\R[x]_d: f|_{\supp\mu}=0\}$.
\end{proposition}

Thus $V_\R(\ker M_d(y))$ always contains the Zariski closure of $\supp\mu$,
and equals it once $d$ is at least the largest degree in a generating set of
the vanishing ideal; this closure lies in $\mathcal Z$, and may be smaller if
$\supp\mu$ misses a component of $\mathcal Z$. We denote the recovered carrier
again by $\mathcal Z$. A numerical basis of $\ker M_d(y)$ mixes generators with
their polynomial multiples, so generators are extracted degree by degree
(\Cref{app:generators}); in our experiments the relations appear in the first
degree in which the kernel is nontrivial, and no reduction is needed. Step~1
does not recover the inequalities \eqref{eq:boundaries}, since
$\ker M_d(y)$ depends only on the Zariski closure of $\supp\mu$.

\subsection{Step 2: sampling the algebraic carrier}
\label{sec:sampling}

\Cref{problem_main} does not permit sampling from $\mu$, but $\mathcal Z$ is
given explicitly by its equations and can be sampled. We use the random-slicing
method of \citet{Breiding2020}, implemented in \texttt{HomotopyContinuation.jl}
\citep{BreidingTimme2018}, which produces points distributed uniformly with
respect to the intrinsic volume of $\mathcal Z$.

Briefly, $\mathcal Z$ is intersected with a random affine subspace
$\mathcal L_{A,b}=\{x:Ax=b\}$ of complementary dimension, where
$A\in\R^{s\times n}$ and $b\in\R^s$; the intersection consists of finitely many
points, which are found by solving one polynomial system. If $(A,b)$ is drawn
from a Gaussian density reweighted by a weighted count of the intersection
points, and one of these points is then chosen with probability proportional to
the reciprocal of its weight, the resulting point is uniformly distributed on
$\mathcal Z$ \citep{Breiding2020}. The weights and the rejection step are
described in \Cref{app:sampling}.

The samples lie exactly on $\mathcal Z$, so no point is accepted merely because
it lies near the support in the ambient space, and since they are uniform in
intrinsic volume, the resolution of the sample depends only on the budget $N$.

\subsection{Step 3: the relative Christoffel--Darboux kernel on the carrier}
\label{sec:cd}

Let $\mu$ and $M_d(y)$ be as above. Under \Cref{problem_main}, $M_d(y)$ is
singular as soon as $d$ is at least the lowest degree of a polynomial vanishing
on $\supp\mu$. The bilinear form
$\langle p,q\rangle_\mu=\int pq\,\de\mu$ on $\R[x]_d$ has nullspace
\begin{equation*}
    I_d:=\ker M_d(y)=\Big\{p\in\R[x]_d:\textstyle\int p^2\,\de\mu=0\Big\},
\end{equation*}
which is the space already considered in \Cref{sec:kernel_method}. The quotient
$\mathcal H_d:=\R[x]_d/I_d$ is a finite-dimensional Hilbert space, which we
identify with a space of polynomial functions on the carrier
$\mathcal Z=V_\R(I_d)$. If $\{P_1,\ldots,P_{r_d}\}$ is an orthonormal basis of
$\mathcal H_d$, the \emph{relative} CD kernel is
$\kappa_d(x,z)=\sum_{j=1}^{r_d}P_j(x)P_j(z)$ for $x,z\in \mathcal Z$, and on
the diagonal
\begin{equation}
\label{eq:CD_kernel}
  \kappa_d(x)=\kappa_d(x,x)=v_d(x)^\top M_d(y)^\dagger v_d(x),
  \qquad x\in \mathcal Z,
\end{equation}
where $M_d(y)^\dagger$ is the Moore--Penrose pseudoinverse; the Christoffel
function is $\Lambda_d=1/\kappa_d$. The quotient construction and the formula
involving the pseudoinverse are developed in
\citet{pauwels2021data,Lasserre_Pauwels_Putinar_2022}. Note that
\eqref{eq:CD_kernel} can be computed from the same data as Step~1, and that it
is well defined, although $M_d(y)$ is singular, because the pseudoinverse is
used.

Given a threshold $\tau_d$, membership is decided by the sublevel set
\begin{equation}
\label{eq:level_set}
    S_{d,\tau_d}:=\left\{x\in \mathcal Z:\kappa_d(x)\leq\tau_d\right\} .
\end{equation}
The optimal threshold $\tau_d^*$ is the one which minimizes
$\dH(S_{d,\tau_d},S)$. When the moment matrix comes from a moment relaxation of
a POP, an estimate $\hat\tau_d$ of $\tau_d^*$ is given by the largest threshold
for which all the retained points $\mathcal M_{\mathrm{adm}}$ are feasible;
we do not use this rule in the experiments.
The choice is motivated by the dichotomy of the Christoffel function:
$\kappa_d$ remains relatively small on the support and grows quickly away from
it, typically exponentially in $d$.

This dichotomy has been proved only under additional assumptions on $S$ and
$\mu$. For full-dimensional supports, \citet{lasserre2019empirical} establish
Hausdorff convergence of suitable sublevel sets of the Christoffel function,
and \citet{vu2022rate} give rates for finite samples. These results require, in
particular, that $S=\overline{\operatorname{int}(S)}$, and so they do not cover
\Cref{problem_main}, in which $S$ has empty interior in $\R^n$. For
lower-dimensional algebraic supports the relative kernel is still well defined
\citep{pauwels2021data}, but in our setting no computable thresholds $\tau_d$
are known which guarantee that $\dH(S_{d,\tau_d},\supp\mu)\to0$. We therefore
use \eqref{eq:level_set} as a heuristic.

\citet{bentancur2026mollified} introduce a mollified CD kernel and prove a
dichotomy of this kind under intrinsic regularity assumptions. Its evaluation,
however, requires a reference measure, mollifiers on the recovered variety and
the associated integrals, none of which is supplied by the moments. We
therefore use the ordinary relative CD kernel, and regard the mollified kernel
as a route to guarantees in future work.

\subsection{Output and choice of parameters}
\label{sec:output}

\Cref{alg:support} returns two objects. The set $S_{d,\tau_d}$ of
\eqref{eq:level_set} is semialgebraic by construction, since $\mathcal Z$ is
algebraic and $\kappa_d$ is a polynomial of degree $2d$. The method therefore
describes the estimated support in the same terms as
\eqref{eq:generators}--\eqref{eq:boundaries}, in addition to providing a finite
set of points. That set, $\mathcal M_{\mathrm{adm}}$, is a sample from
$S_{d,\tau_d}$, and the Hausdorff distances reported in
\Cref{sec:experiments} are computed for it.

Two parameters are not determined by the theory. There is no formula for
$\tau_d$ in terms of the moments, and we choose it from the distribution of the
values of $\kappa_d$ on $\mathcal M$, which is bimodal whenever the dichotomy
is pronounced (\Cref{fig:cone}, left). Nor can the smallest $d$ achieving a
given $\delta$ be predicted, and in any case the order of the relaxation bounds
$d$ from above. \Cref{alg:support} therefore does not guarantee
\eqref{eq:goal} for a prescribed $\delta$; the experiments report the value of
$\delta$ attained.

\section{Experiments}
\label{sec:experiments}

We test \Cref{alg:support} in two settings: with moments computed in closed
form for measures which satisfy \Cref{problem_main} exactly
(\Cref{sec:analytic}), and with pseudo-moment sequences obtained from POP
relaxations (\Cref{sec:pop}). In the second setting the hypotheses need not
hold: the sequence need not be a moment sequence, and if it is one, the
representing measure may be supported on any subset of the set of minimizers.
We therefore regard the POP experiments as exploratory; they test whether the
method extracts useful geometric information beyond the setting of exact
moments.

\subsection{Setup and baseline}
\label{sec:setup}

In all experiments the support is estimated by thresholding the relative CD
kernel \eqref{eq:CD_kernel}; the two settings differ only in the origin of the
moments. The comparison concerns the points at which the kernel is evaluated.
Our \textbf{carrier} estimator evaluates $\kappa_{d_v}$ on
$\mathcal M\subset \mathcal Z$, sampled as in \Cref{sec:sampling}. The
\textbf{grid} estimator, which serves as the baseline, evaluates $\kappa_{d_g}$
on a uniform grid in the ambient box \eqref{eq:boundaries_given}; it is CD
thresholding without Steps~1 and~2. Both estimators use moments of the same
measure or relaxation, and in every case the grid is allowed the higher degree
$d_g\ge d_v$. Off $\mathcal Z$ the pseudoinverse ignores the component of
$v_d(x)$ in $\ker M_d(y)$; a regularized kernel built from
$(M_d(y)+\varepsilon I)^{-1}$ would penalize that component, and is a natural
further baseline which we have not evaluated. In each example $\tau_d$ is
chosen by hand, usually as a percentile of the values of $\kappa_d$ on the
samples from the variety, and the reported errors depend on this choice. For each estimator we record the degree, the
number of samples and the Hausdorff distance $\dH$ between the accepted points
and the true support, computed against a fine discretization of $S$. The
results are collected in \Cref{tab:recovery-experiments}.

\begin{table}[t]
\centering
\small
\setlength{\tabcolsep}{5pt}
\begin{tabular}{>{\raggedright\arraybackslash}p{4.3cm}ccc c ccc}
\toprule
& \multicolumn{3}{c}{\textbf{\Cref{alg:support}}} & & \multicolumn{3}{c}{\textbf{CD, ambient grid}} \\
\cmidrule(lr){2-4}\cmidrule(lr){6-8}
Experiment & $d_v$ & Samples & $\dH$ & & $d_g$ & Samples & $\dH$ \\
\midrule
\multicolumn{8}{l}{\emph{Moments in closed form} (\Cref{sec:analytic})}\\[2pt]
\makecell[l]{Parabola\\ $y=x^2,\ 0\leq x\leq1$}
& 2 & $1.1\times10^4$ & \textbf{0.014} & & 4 & $10^6$ & 0.095 \\[4pt]
\makecell[l]{Hyperbola\\ $xy=1,\ 0.5\leq x\leq1$}
& 2 & 4945 & \textbf{0.002} & & 4 & $10^6$ & 0.114 \\[4pt]
\makecell[l]{Hyperbola 2\\ $x^2-y^2=1,\ 1\leq x\leq2$}
& 2 & $1.4\times 10^4$ & \textbf{0.003} & & 4 & $10^6$ & 0.096 \\[4pt]
\makecell[l]{3D cone (Section \ref{sec:analytic})\\ $z^2-x^2-y^2=0$,\\ $z\in[-1,-0.5]\cup[0.5,1]$}
& 6 & $2\times 10^4$ & \textbf{0.090} & & 8 & $10^6$ & 0.240 \\[4pt]
\makecell[l]{Cylinder in 3D\\ $x^2+y^2=1,\ z\in[0,1]$}
& 2 & $2.3\times10^4$ & \textbf{0.058} & & 8 & $10^6$ & 0.982 \\
\midrule
\multicolumn{8}{l}{\emph{Pseudo-moments from POP relaxations} (\Cref{sec:pop})}\\[2pt]
\makecell[l]{Line segment (Section \ref{app:segment})\\ $y=x,\ -0.5\leq x\leq0.5$}
& 1 & $1.1\times10^4$ & \textbf{0.002} & & 2 & $10^6$ & 0.003 \\[4pt]
\makecell[l]{Semicircle (Section \ref{app:semicircle})\\ $x^2+y^2=1,\ y\geq0$}
& 2 & $2.2\times10^4$ & \textbf{0.050} & & 4 & $10^6$ & 0.270 \\
\makecell[l]{Product of segments in 10D \\ (\Cref{sec:pop-product-segments})}
& 1 & $5\times 10^4$ & \textbf{0.530}& & 2& $9.7\times10^6$ & $1.32$   \\
\makecell[l]{Product measure in 4D\\ (\Cref{app:mixed})}
&2 & $10^4$ & \textbf{0.1} & &4& $6.2\times10^6$ & 0.246 \\
\bottomrule
\end{tabular}
\caption{Support recovery by \Cref{alg:support} and by CD thresholding on an
ambient grid. Here $d_v$ and $d_g$ are the kernel degrees, \emph{Samples} is the
number of points at which the kernel is evaluated, and $\dH$ is the Hausdorff
distance to the true support, computed against a fine discretization of
$\supp\mu$. The grid baseline uses the higher degree and more samples in every row, and has the
larger error in every row.}
\label{tab:recovery-experiments}
\end{table}

\subsection{Cone moments computed in closed form}
\label{sec:analytic}

\begin{figure}[t]
    \centering
    \includegraphics[width=0.46\linewidth]{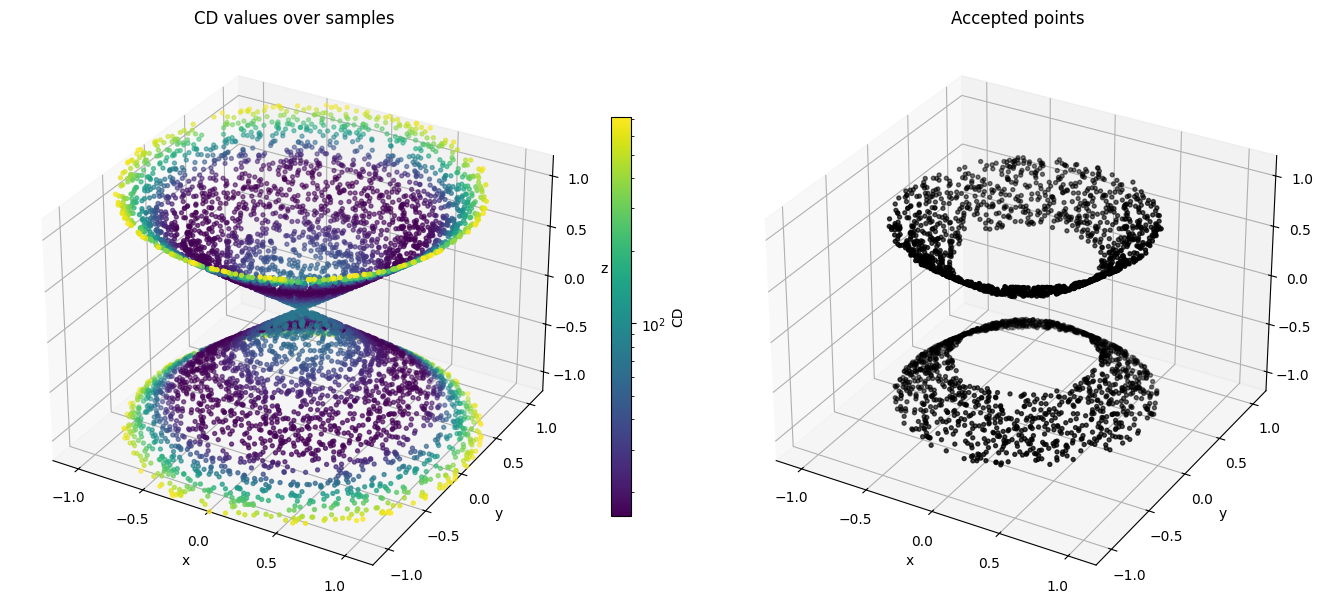}
    \caption{The cone, with moments in closed form. \emph{Left:} values of
    $\kappa_6$ at $20000$ points sampled uniformly on the recovered carrier
    $\{x_3^2-x_2^2-x_1^2=0\}$, on a logarithmic color scale; the dichotomy
    appears as dark inner annuli and bright rims. \emph{Right:} the $10600$
    points accepted at $\tau_6=82.7$; both connected components are
    recovered.}
    \label{fig:cone}
\end{figure}

Here $\mu$ is the normalized uniform measure, with respect to intrinsic surface
area, on a set which satisfies \Cref{problem_main} exactly, and its moments are
computed in closed form. We describe the cone, which is the most intricate of
the five examples: it has codimension one in $\R^3$, it has a singular point at
the origin, and its support, which is the point of the experiment, has
\emph{two} connected components,
\begin{equation*}
    S=\left\{(x_1,x_2,x_3)\in\R^3:\ x_3^2-x_2^2-x_1^2=0,\
      x_3\in[-1,-0.5]\cup[0.5,1]\right\} .
\end{equation*}
Step~1 finds the polynomial $f_1(x)=x_3^2-x_2^2-x_1^2$ in $\ker M_2(y)$; this
is the polynomial which defines $\mathcal Z$ in \eqref{eq:generators}. The singular point
of the cone lies outside $S$ and is met by random slices with probability zero.
The CD kernel of degree $6$, which requires $M_6(y)$, is then evaluated at
$2\times10^4$ samples on the cone and thresholded at $\tau_6=82.7$, and $10600$
points are retained. The projection of the accepted set on the $x_3$-axis is
$[-0.958,-0.489]\cup[0.489,0.958]$, against the true $[-1,-0.5]\cup[0.5,1]$.
Both components are recovered, and recovered separately, although the method is
not told how many there are; the error is concentrated at the outer rims, where
the density is lowest (\Cref{fig:cone}). The Hausdorff distance
$\dH(S^{\mathrm{variety}}_{6,\tau_6},S)=0.09$ was computed against an
$\varepsilon$-net of $S$ constructed by hand, with $\varepsilon=0.01$.
Thresholding $10^6$ grid points gives
$\dH(S^{\mathrm{grid}}_{8,\tau_8},S)=0.240$.

\subsection{Pseudo-moments from polynomial optimization}
\label{sec:pop}

The set of minimizers of a POP is semialgebraic, and when it is infinite it has
positive dimension; this is the setting of \Cref{problem_main}. We solve the
POPs by the moment relaxation \citep{Lasserre2001GlobalOptimization}, as
implemented in TSSOS \citep{magron2021tssos,wang2021tssos}. Under the usual
assumptions of the moment--SOS hierarchy the relaxations converge to the global
optimal value, but they do not return the minimizers: one obtains a finite
pseudo-moment sequence, from which the minimizers must be located. When the
objective and constraints are known, samples on $\mathcal Z$ could also be
filtered by them directly; we do not do so, since our aim is to test what can
be recovered from the pseudo-moments alone. When the flat extension property
\citep{CurtoFialkow1998} holds, the GNS construction
\citep{LopezQuijorna2021DetectingGNS} solves the extraction problem. In the
examples below the flat extension property fails, and the GNS construction
cannot be used.

\subsubsection{The polynomial optimization setting}
\label{sec:POP_setting}

In each of the following examples we solve a polynomial optimization problem
\begin{align*}
    \min_{x\in\R^{n}} \quad & F(x) \\
    \text{s.t.}\quad & g_j(x) \geq 0, \quad j=1,2,\ldots,l,
\end{align*}
whose set of minimizers $S$ is known in closed form. The relaxation returns a
moment matrix $M_d(y)$, which we pass to \Cref{alg:support} as though it had
been generated by a measure supported on $S$. With $F(x)=(x_1-x_2)^2$ and
suitable constraints $g_j$, the set $S$ is a line segment
(\Cref{app:segment}); with $F(x)=(x_1^2+x_2^2-1)^2$ it is a subset of the
circle (\Cref{app:semicircle}). In both cases $S$ is recovered to the accuracy
reported in \Cref{tab:recovery-experiments}.

\subsubsection{A product of segments in $2n$ dimensions}
\label{sec:pop-product-segments}

\begin{figure}[t]
    \centering
    \includegraphics[width=0.6\linewidth]{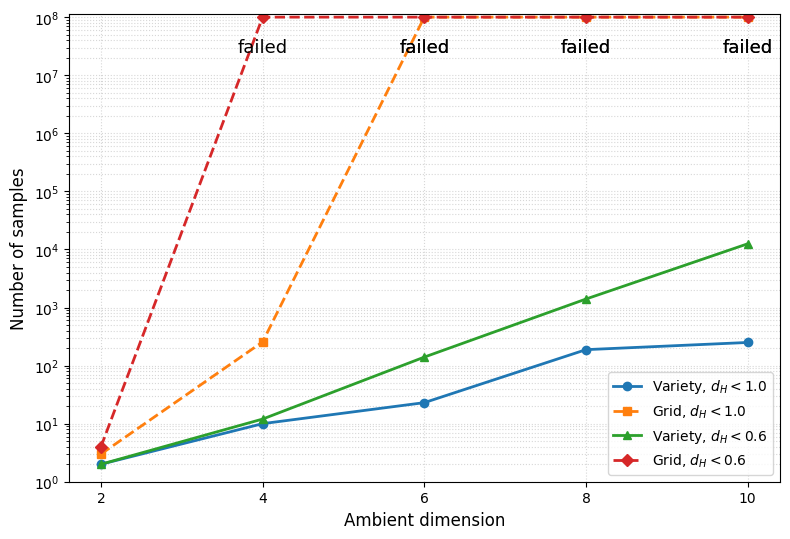}
    \caption{Number of evaluation points needed to reach Hausdorff errors
    $\dH<1.0$ and $\dH<0.6$. The carrier estimator uses $\kappa_1$ on
    $\mathcal{Z}$; the grid baseline is allowed $\kappa_3$ in $\R^{2n}$.
    Failed runs are plotted at $10^8$ for display only; they indicate failure
    within the tested budget and are not estimates of the number of samples
    needed.}
    \label{fig:10d}
\end{figure}

To examine how the approach scales, we consider a family of examples in
dimensions $2n\in\{2,4,6,8,10\}$. We solve, at relaxation degree $2$,
\begin{align}
\label{eq:pop_2nd}
    \min_{x\in\R^{2n}} \quad & F(x)=\sum_{i=1}^{n}(x_{2i-1}-x_{2i})^2 \\ \nonumber
    \text{s.t.}\quad & -0.5\le x_{j}\le0.5,\qquad j=1,\ldots,2n,
\end{align}
whose set of minimizers is $S^n\subset\R^{2n}$, where
$S=\{(x_1,x_2)\in\R^2: x_1=x_2,\ -0.5\le x_1\le0.5\}$. The objective is chosen
to be quadratic, so that its evaluation requires moments only up to degree $2$.
A dense moment matrix $M_d(y)$ of order $d$ has $\binom{d+2n}{d}^2$ entries, and
raising the relaxation degree soon becomes prohibitively expensive. We study
how the number of evaluation points needed to reach a prescribed Hausdorff
accuracy grows with the dimension. The Hausdorff distances in \Cref{fig:10d}
are computed against a tensor-product discretization of $S^n$ with $20$ points
per intrinsic coordinate.

For the stricter target $\dH<0.6$, the number of carrier samples required grows
by about one order of magnitude for each additional intrinsic dimension over
the range tested. For the grid baseline we also solve the POP at relaxation
degree $3$ and allow it to use $\kappa_3$. It succeeds for $2n=2$; for $2n=4$
it needs $256$ grid points to reach $\dH\le1$, and with up to $10^8$ grid points
it does not reach $\dH<0.6$ for any of the thresholds tried. For $2n\ge8$ memory is the limit: a grid with $9$ points per coordinate
has $9^8\approx4.3\times10^7$ points, and evaluating the degree-$2$ CD kernel
on it needs about $115$\,GiB. The CD kernel of degree $1$ cannot be used
on the grid, because it always accepts the point $x=(-1,1,-1,1,\ldots)$ (see
\Cref{fig:segment_deg_1}), and so forces the Hausdorff distance to be at least
$\dist(x,S^n)=\sqrt{2n}$. Restriction to the carrier thus reduces the
dimension to be discretized from $2n$ to $n$.

A four-variable example with quartic terms, the product of a semicircle and a
segment, is given in \Cref{app:mixed}.\label{sec:pop-product-mixed}

\subsection{The effect of restricting to the carrier}
\label{sec:ablation}

\begin{figure}[t]
    \centering
    \includegraphics[width=0.82\linewidth]{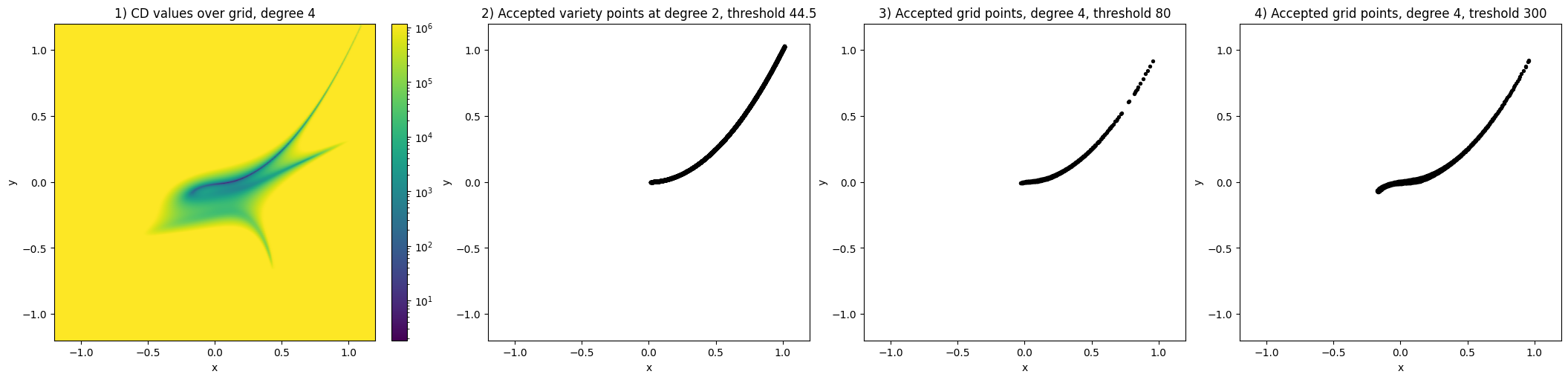}
    \caption{The parabola $y=x^2$, $0\le x\le 1$, with moments in closed form.
    \emph{1)} The CD kernel of degree $4$ on the grid. \emph{2)} Points of
    $\mathcal Z$ accepted by \Cref{alg:support}, $\dH=0.014$. \emph{3)} Grid
    points accepted at threshold $80$; the estimate has large gaps,
    $\dH=0.095$. \emph{4)} Grid points accepted at the higher threshold $500$;
    the gaps are filled, but points far from the support are admitted,
    $\dH=0.17$.}
    \label{fig:parabola}
\end{figure}

\Cref{tab:recovery-experiments} shows that the carrier estimator has a smaller
Hausdorff error than the grid in all nine experiments, and
\Cref{fig:parabola} illustrates the reason. Samples drawn on $\mathcal Z$
satisfy its defining equations to machine precision, with residuals of order
$10^{-15}$, whereas even grid points close to $\mathcal Z$ have residuals of
order $10^{-3}$ to $10^{-1}$. In our experiments the region in which the CD
kernel is small forms a thin neighbourhood of the support, which a grid may
fail to resolve. A small threshold then leaves gaps in the estimated support; a
larger threshold fills the gaps, but admits points farther from $\mathcal Z$
and so increases the Hausdorff error. Evaluating the kernel only on
$\mathcal Z$ avoids this compromise.

In dimension ten (\Cref{sec:pop-product-segments}) the resolution of the grid
itself becomes a limitation. Indeed, the Hausdorff error of any estimator which selects points
from a grid can be bounded below, whatever function is used to make the
selection. Let
\[
G_m=\{-1,-1+\tfrac{2}{m},\ldots,1\}^n
\]
denote the uniform grid of $(m+1)^n$ points in $[-1,1]^n$.

\begin{definition}[Dispersion]
\label{def:disp}
Let $G\subset\R^n$ be finite and $S\subseteq\R^n$ compact. The
\emph{dispersion} of $G$ relative to $S$ is
\begin{equation*}
  \disp(G,S):=\sup_{y\in S}\ \dist(y,G)=\sup_{y\in S}\ \min_{x\in G}\norm{x-y},
\end{equation*}
that is, the smallest $r$ for which $G$ is an $r$-net of $S$
\citep[\S4.2]{Vershynin2018}; it is also known as the \emph{fill distance}
\citep{Niederreiter1992}. It is the one-sided Hausdorff distance from $S$ to
$G$, and it does not depend on how points of $G$ are subsequently selected.
\end{definition}

\begin{theorem}
\label{thm:grid}
Let $S\subseteq \mathcal Z\cap[-1,1]^n$ be compact.
\begin{enumerate}[label=\textup{(\roman*)}]
\item For every $m$ and every $\hat S\subseteq G_m$,
  $\dH(\hat S,S)\ge\disp(G_m,S)$.
\item Let $\mathcal M\subset \mathcal Z$ be finite and
  $\hat S=\{x\in\mathcal M:\dist(x,S)\le\rho\}$. If $\rho\ge\disp(\mathcal M,S)$,
  then $\dH(\hat S,S)\le\rho$.
\item Let $S_n=\{t\mathbf 1:|t|\le\tfrac12\}\subset\R^n$, the set of minimizers
  of $\sum_{i<j}(x_i-x_j)^2$ over $[-\tfrac12,\tfrac12]^n$, with carrier
  $\mathcal Z=\R\mathbf 1$. Then $\disp(G_m,S_n)\ge\sqrt n/m$ for every $m\ge1$.
  On the other hand, if $\mathcal M$ consists of $N\ge4\ln(N/\epsilon)$ points
  drawn uniformly from $\mathcal Z\cap[-1,1]^n$ and $\hat S=\mathcal M\cap S_n$,
  then with probability at least $1-\epsilon$
  \[
    \dH(\hat S,S_n)\ \le\ 4\sqrt n\,\frac{\ln(N/\epsilon)}{N}.
  \]
\end{enumerate}
\end{theorem}

\Cref{thm:grid} concerns only the choice of evaluation points. Part~(i) holds
for every selection rule, and in particular for CD thresholding; parts~(ii)
and~(iii) use the ideal rule $\dist(x,S)\le\rho$, which is not available in
practice, so they bound what the carrier can achieve, not what
\Cref{alg:support} achieves. Whether CD thresholding comes close to the ideal
rule is the question tested by the experiments; $S_n$ is a simplified model of
\Cref{sec:pop-product-segments}. Thus, for the same budget of $N=(m+1)^n$
evaluations, the two errors on $S_n$
satisfy $\dH^{\mathrm{grid}}/\dH^{\mathrm{carrier}}\ge
(m+1)^n\big/\big(4m\ln((m+1)^n/\epsilon)\big)$\label{eq:ratio}, a ratio which
grows exponentially in $n$ for every fixed $m$. To bring the error below a
prescribed $\delta<\sqrt n=\operatorname{diam}S_n$, the grid needs at least
$(1+\sqrt n/\delta)^n$ evaluations, and the carrier needs
$O\big(\sqrt n\,\delta^{-1}\ln(\sqrt n/\delta\epsilon)\big)$. The count is of evaluations, not of computing time: each sample on
$\mathcal Z$ costs one polynomial-system solve. The proof is in
\Cref{app:grid}; part (iii) does not depend on the grid being cubic
(\Cref{rem:anydesign}).

\section{Conclusion and limitations}
\label{sec:conclusion}
\label{sec:limitations}

We have given a method for estimating, from its moments alone, a support which
is compact, semialgebraic and of positive codimension. The kernel of the moment
matrix recovers the equations, and the relative Christoffel--Darboux kernel,
evaluated on a uniform sample of the recovered carrier, recovers the
inequalities. In our experiments CD thresholding on an ambient grid has the
larger Hausdorff error in every case, although it is run at a higher degree;
the margin is negligible for the segment and large elsewhere. For a model
problem, and with an ideal selection rule on the carrier, \Cref{thm:grid} shows
that the gap in attainable error grows exponentially with the dimension.

Several questions remain open. Step~1 is exact for true moments, but for the
pseudo-moments of \Cref{sec:pop} the kernel of $M_d(y)$ is computed
numerically and its singular-value cut-off is set by inspection; a
perturbation theory along the lines of
\citet{KlepPovhVolcic2018_MinimizerExtractionIsRobust} is needed. Step~3 relies
on a dichotomy proved for full-dimensional supports
\citep{lasserre2019empirical,vu2022rate} and, for the mollified kernel, on
varieties \citep{bentancur2026mollified}; for the ordinary relative kernel in
codimension at least one we know of no computable thresholds for which
$\dH(S_{d,\tau_d},\supp\mu)\to0$, hence no rate and no principled choice of
degree or threshold, and we have not studied the sensitivity of the results to
$\tau_d$. The sampler assumes a smooth carrier, and the pseudoinverse of order
$\sigma(d)=\tbinom{n+d}{d}$ dominates the cost.

\bibliography{library}
\bibliographystyle{plainnat}

\appendix

\section{From kernel elements to generators}
\label{app:generators}

\begin{remark}
\label{rem:generators}
Even when the vanishing ideal $I(\supp\mu)$ is generated by a single polynomial
$f$ of degree $d$, the polynomial $f$ can be read off directly from the kernel
of the moment matrix only in the first degree in which a nontrivial relation
appears. Indeed, if $I(\supp\mu)\cap\R[x]_{d-1}=\{0\}$ and
$I(\supp\mu)\cap\R[x]_d=\operatorname{span}\{f\}$, then $\dim\ker M_d(y)=1$, and
$f$ may be recovered, up to a scalar factor, from a basis vector of
$\ker M_d(y)$. In higher degrees the kernel also contains polynomial multiples
of relations already found: $f,x_1f,\ldots,x_nf\in\ker M_{d+1}(y)$, so that
$\dim\ker M_{d+1}(y)>1$ although the vanishing ideal may still have a single
generator. More generally, a numerical basis of $\ker M_d(y)$ consists of
arbitrary linear combinations of generators of the vanishing ideal and of their
polynomial multiples.

To extract generators one proceeds degree by degree, separating the new
relations in each degree from those generated in lower degrees. Equivalently,
in degree $d$ one removes from $\ker M_d(y)$ the subspace generated by
polynomial multiples of kernel elements of lower degree, and the remaining
independent relations give new generators. In computational algebra this is
done by row reduction with respect to a monomial ordering, by Gr\"obner-basis
techniques, or by border-basis methods; see, for instance,
\citet{laurent2009generalized,lasserre2013moment}.
\end{remark}

\begin{remark}[Several generators]
\label{rem:several}
The same reasoning applies, under suitable assumptions, when $\mathcal Z$ is
defined by several polynomials. Suppose that the first nontrivial relations in
$I(\supp\mu)$ occur in degree $d$, and that
\[
    I(\supp\mu)\cap\mathbb R[x]_d
    =
    \operatorname{span}\{f_1,\ldots,f_k\}.
\]
Then $\dim\ker M_d(y)=k$. A numerical nullspace computation may return a
different basis $v_j=\sum_{i=1}^k C_{ji}f_i$ with $C$ invertible, but then
$\{v_1=\cdots=v_k=0\}=\{f_1=\cdots=f_k=0\}$, and the same algebraic set is
obtained. In polynomial optimization the $f_i$ often involve disjoint cliques of
variables, and each independent regular equation then typically lowers the
local dimension by one. Since samples are drawn only on the recovered carrier,
each such equation removes, before CD thresholding, one ambient dimension in
which spurious samples could otherwise lie.
\end{remark}

\section{The random-slicing sampler}
\label{app:sampling}

We recall the construction of \citet{Breiding2020} used in
\Cref{sec:sampling}. Suppose that $\mathcal Z$ is a smooth algebraic manifold of
dimension $s$, defined by $F=(f_1,\ldots,f_k)$ as in \eqref{eq:generators}, and
of positive, finite volume with respect to the induced volume form. For
$A\in\R^{s\times n}$ and $b\in\R^s$ put $\mathcal L_{A,b}=\{x\in\R^n:Ax=b\}$,
and define
\begin{equation}
\label{eq:slicing-count}
  \bar 1(A,b)=\sum_{x\in \mathcal Z\cap\mathcal L_{A,b}}\frac{1}{\alpha(x)},
  \qquad \bar 1:\R^{s\times n}\times\R^{s}\to\R ,
\end{equation}
where $\alpha(x)>0$ is a normalizing factor computed from the Jacobian of $F$ at
$x$. If $\varphi$ denotes the standard Gaussian density on
$\R^{s\times n}\times\R^s$, the integral-geometric identity
$\mathbb E_\varphi[\bar 1]=\vol(\mathcal Z)$ holds. The sampler draws $(A,b)$
from the density
\begin{equation*}
  \psi(A,b)=\frac{\varphi(A,b)}{\vol(\mathcal Z)}
            \sum_{x\in \mathcal Z\cap\mathcal L_{A,b}}\alpha(x)^{-1}
\end{equation*}
by rejection: one proposes $(A,b)\sim\varphi$ and accepts it with probability
$\eta\,\bar 1(A,b)$, where the constant $\eta$ depends only on $\deg\mathcal Z$
and on a bound for $\norm{x}$ over $\mathcal Z$ \citep{Breiding2020}. Each
evaluation of $\bar 1$ requires the solution of one polynomial system. Given an
accepted $(A,b)$, one returns $x\in \mathcal Z\cap\mathcal L_{A,b}$ with
probability $\alpha(x)^{-1}/\bar 1(A,b)$, and the resulting point is uniformly
distributed on $\mathcal Z$.

\section{Proof of \Cref{thm:grid}}
\label{app:grid}

Throughout, $\dist(y,G)=\min_{x\in G}\norm{x-y}$, and $\dH(\emptyset,S)=+\infty$.

\emph{(i).} For $\hat S\subseteq G_m$ and $y\in S$ we have
$\dist(y,\hat S)\ge\dist(y,G_m)$. Hence
$\dH(\hat S,S)\ge\sup_{y\in S}\dist(y,\hat S)\ge\disp(G_m,S)$. In particular,
this applies to $\hat S=G_m\cap\{\kappa_d\le\tau\}$ for every $d$ and $\tau$.

\emph{(ii).} Let $y\in S$. By hypothesis there is $x\in\mathcal M$ with
$\norm{x-y}\le\disp(\mathcal M,S)\le\rho$; then $\dist(x,S)\le\rho$, so that
$x\in\hat S$ and $\dist(y,\hat S)\le\rho$. Conversely, every $x\in\hat S$
satisfies $\dist(x,S)\le\rho$ by construction. Both one-sided terms of $\dH$
are therefore at most $\rho$.

\emph{(iii): recovery of the carrier.} An affine form $c_0+\sum_ic_ix_i$
vanishes on $S_n$ if and only if $c_0=0$ and $\sum_ic_i=0$, so the affine part
of $\ker M_1(y)$ is spanned by $x_i-x_1$, $i=2,\ldots,n$, whose common zero set
is the line $\R\mathbf 1$. Step~1 therefore returns $\mathcal Z$ exactly in
degree one, as in the case $n=2$ of \Cref{app:segment}.

\emph{(iii): the grid.} The coordinates of points of $G_m$ are $-1+2j/m$,
$j=0,\ldots,m$, and the centres of its cells have coordinates $-1+(2j+1)/m$,
$j=0,\ldots,m-1$. Put $\gamma=0$ if $m$ is odd and $\gamma=1/m$ if $m$ is even.
In either case $\gamma$ is a coordinate of a cell centre and
$|\gamma|\le\tfrac12$, so that $c=\gamma\mathbf 1\in S_n$. Each coordinate of
$c$ is at distance exactly $1/m$ from the nearest grid coordinate, whence
$\dist(c,G_m)=\sqrt n/m$ and $\disp(G_m,S_n)\ge\sqrt n/m$. Relative to
$\operatorname{diam}S_n=\sqrt n$, the error of the grid is thus at least $1/m$,
that is, $1/(N^{1/n}-1)$ with $N=(m+1)^n$ evaluations; it tends to zero only if
the budget grows faster than every exponential in $n$.

\emph{(iii): the carrier.} We have $\mathcal Z\cap[-1,1]^n=\{t\mathbf 1:|t|\le1\}$,
and a uniform point of this set is $t\mathbf 1$ with $t$ uniform on $[-1,1]$.
Partition $[-\tfrac12,\tfrac12]$ into $K$ intervals of length $1/K$. A sample
falls in a given interval with probability $1/(2K)$, so the probability that
some interval receives none of the $N$ samples is at most
$K(1-\tfrac1{2K})^N\le Ke^{-N/(2K)}$. On the complementary event every
$t\in[-\tfrac12,\tfrac12]$ lies within $1/K$ of some sample
$t_i\in[-\tfrac12,\tfrac12]$; the point $t_i\mathbf 1$ belongs to $S_n$, and
hence to $\hat S$, and $\dist(t\mathbf 1,\hat S)\le\sqrt n/K$. Since
$\hat S\subseteq S_n$, the other one-sided term vanishes. Hence
$\dH(\hat S,S_n)\le\sqrt n/K$ on that event. Take
$K=\lfloor N/(2\ln(N/\epsilon))\rfloor$. Then $e^{-N/(2K)}\le\epsilon/N$, so
the exceptional probability is at most $\epsilon K/N\le\epsilon$; and
$N\ge4\ln(N/\epsilon)$ gives $K\ge N/(4\ln(N/\epsilon))$, whence
$\sqrt n/K\le4\sqrt n\ln(N/\epsilon)/N$.

\emph{Consequences.} At $N=(m+1)^n$ the quotient of the two bounds is
$(\sqrt n/m)\big/\big(4\sqrt n\ln(N/\epsilon)/N\big)=N/(4m\ln(N/\epsilon))$,
which is the ratio stated after \Cref{thm:grid}. For a target $\delta<\sqrt n$,
(i) and (iii) require $\sqrt n/m\le\delta$, that is, $m\ge\sqrt n/\delta>1$
and $N=(m+1)^n\ge(1+\sqrt n/\delta)^n$. For the carrier, put
$A=4\sqrt n/\delta$ and $N=\lceil2A\ln(2A/\epsilon)\rceil$; since
$\ln\ln(2A/\epsilon)\le\ln(2A/\epsilon)$, we have $A\ln(N/\epsilon)\le N$, and
so the bound in (iii) is at most $\delta$.

\begin{remark}[Arbitrary ambient designs]
\label{rem:anydesign}
The lower bound does not depend on the grid being cubic. If
$G\subset[-1,1]^n$ is any set of $N$ points and $\omega_n$ is the volume of the
unit ball, then the $N$ balls of radius $r$ centred at the points of $G$ cannot
cover the box unless $N\omega_nr^n\ge2^n$. Hence some $c\in[-1,1]^n$ satisfies
$\dist(c,G)\ge2(N\omega_n)^{-1/n}$, and $\omega_n^{-1/n}\sim\sqrt{n/(2\pi e)}$.
Any segment $S$ through $c$, with $\mathcal Z$ the line containing it, then
satisfies $\disp(G,S)\ge2(N\omega_n)^{-1/n}$, and (i) applies. Thus any fixed
set of $N$ evaluation points in the ambient space has dispersion of order
$N^{-1/n}$, whereas $N$ uniform points on $\mathcal Z$ have dispersion of order
$N^{-1/s}$, up to a logarithmic factor (\Cref{rem:generalZ}).
\end{remark}

\begin{remark}[General carriers]
\label{rem:generalZ}
The bound for the carrier in (iii) is the one-dimensional case of a standard
covering argument. Suppose that $\mathcal H^s(\mathcal Z\cap B(y,r))\ge c\,r^s$
for all $y\in S$ and $r\le r_0$, and let $S$ be covered by $\mathcal N$ balls of
radius $\delta/2\le r_0$ with centres in $S$. Then $N$ uniform samples from
$\mathcal Z\cap[-1,1]^n$, a set of $\mathcal H^s$-measure $V$, satisfy
$\disp(\mathcal M,S)\le\delta$ with probability at least $1-\epsilon$ as soon as
$N\ge\frac{V}{c(\delta/2)^s}\ln\frac{\mathcal N}{\epsilon}$
\citep[cf.][]{NiyogiSmaleWeinberger2008}, and (ii) then bounds $\dH$ by
$\max(\delta,\rho)$. The budget grows like $\delta^{-s}$, compared with
$\delta^{-n}$ for the grid.
\end{remark}

\section{The line segment}
\label{app:segment}

The polynomial optimization problem for the line segment of
\Cref{sec:POP_setting}, which we solve at relaxation order~$2$, is
\begin{align}
\label{eq:pop_segment}
    \min_{x_1,x_2\in\R} \quad & F(x_1,x_2)=(x_1-x_2)^2 \\ \nonumber
    \text{s.t.}\quad & x_2+0.5\ge0,\quad 0.5-x_2\ge0,\quad
                       x_1+0.5\ge0,\quad 0.5-x_1\ge0 ,
\end{align}
and its set of minimizers is
$S=\{(x_1,x_2)\in\R^2: x_1=x_2,\ -0.5\le x_1\le0.5\}$. The kernel of $M_1(y)$
is spanned by the vector $(0,\,0.70710678,\,-0.70710678)$, which gives, up to a
scalar factor, the polynomial $f_1(x)=x_2-x_1$ of \eqref{eq:generators}
defining the carrier $\mathcal{Z}$. The line is sampled over the box
$[-1.2,1.2]^2$ with \texttt{HomotopyContinuation.jl} \citep{BreidingTimme2018},
and the grid consists of $1000\times1000$ points in $[-1,1]^2$. We now
examine the behaviour of the CD kernel as the degree $d$ increases.

\begin{figure}[h]
    \centering
    \includegraphics[width=0.8\linewidth]{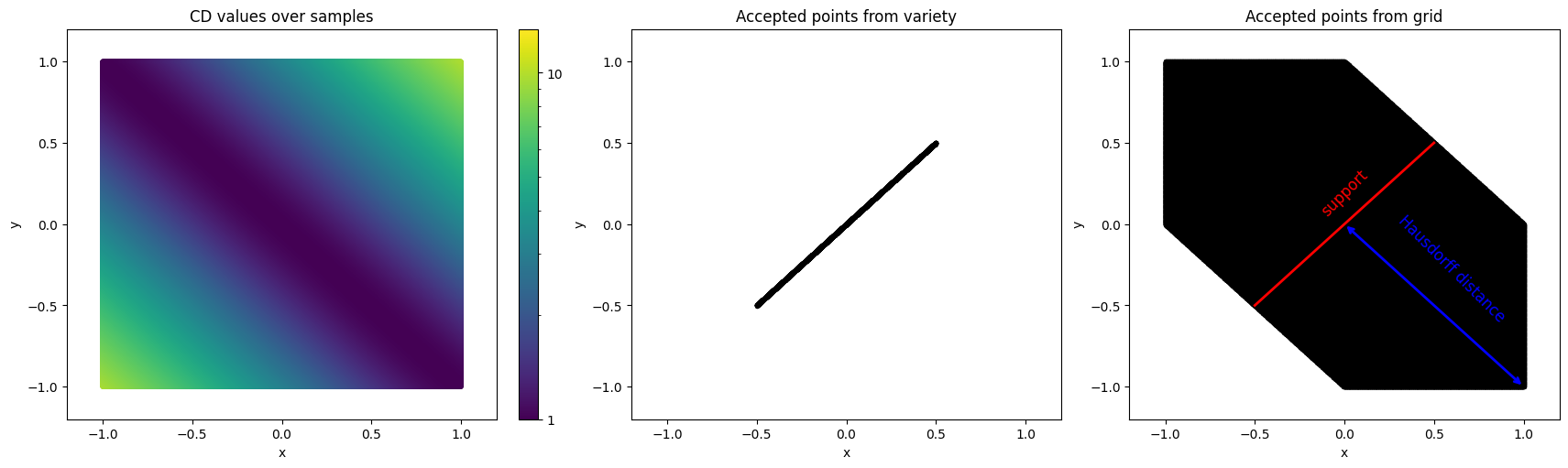}
    \caption{Thresholding the CD kernel $\kappa_1$ on the variety and on the
    grid. \emph{Left:} values of $\kappa_1$ on the $1000\times1000$ grid.
    \emph{Middle:} points of the variety accepted at $\tau_1=3.36$.
    \emph{Right:} grid points accepted at $\tau_1=3.36$.}
    \label{fig:segment_deg_1}
\end{figure}

\emph{Degree one.} The CD kernel $\kappa_1$ of \eqref{eq:CD_kernel} is
evaluated at $1.1\times 10^4$ samples on $\mathcal{Z}$ and at $10^6$ grid
points. With the threshold $\tau_1=3.36$, $4873$ points
$S^{\mathrm{variety}}_{1,\tau_1}$ of $\mathcal Z$ and $659528$ grid points
$S^{\mathrm{grid}}_{1,\tau_1}$ are accepted, and
$\dH(S^{\mathrm{variety}}_{1,\tau_1},S)=0.002$, whereas
$\dH(S^{\mathrm{grid}}_{1,\tau_1},S)=\sqrt{2}$. The kernel $\kappa_1$ and the
accepted points are shown in \Cref{fig:segment_deg_1}.

\begin{figure}[h]
    \centering
    \includegraphics[width=0.95\linewidth]{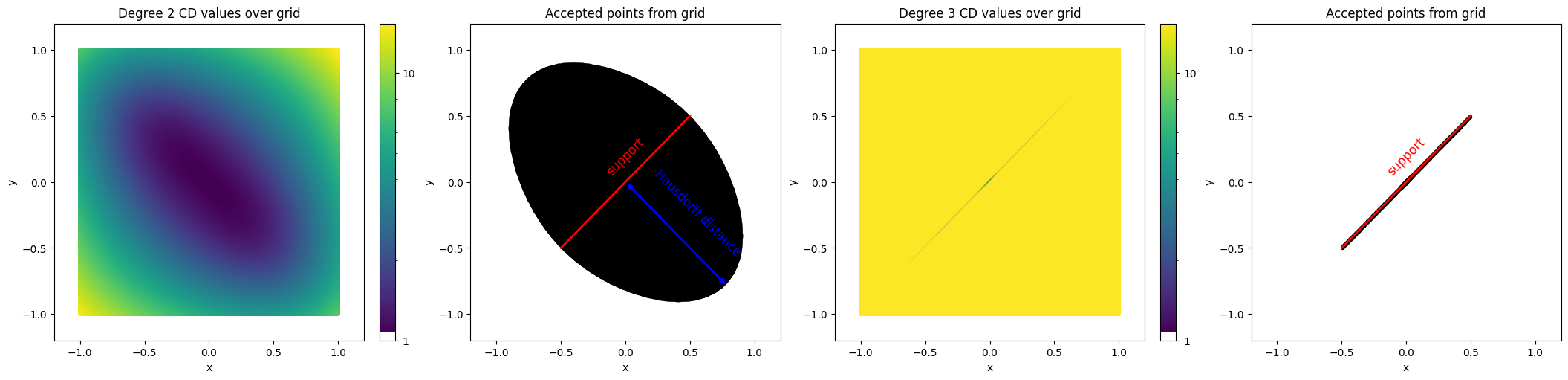}
    \caption{The CD kernels $\kappa_2$ and $\kappa_3$, evaluated and
    thresholded on the $1000\times1000$ grid. \emph{Left:} values of
    $\kappa_2$. \emph{Middle left:} grid points accepted at $\tau_2=3.36$.
    \emph{Middle right:} values of $\kappa_3$. \emph{Right:} grid points
    accepted at $\tau_3=5.0$.}
    \label{fig:segment_deg_2}
\end{figure}

\emph{Degrees two and three.} Since $S^{\mathrm{variety}}_{1,\tau_1}$ already
approximates $S$ sufficiently well, we evaluate the kernels $\kappa_2$ and
$\kappa_3$ on the grid only. The computation of $\kappa_3$ requires $M_3(y)$,
and so \eqref{eq:pop_segment} must be solved again at relaxation order $3$. The
Hausdorff distances are $\dH(S^{\mathrm{grid}}_{2,\tau_2},S)=1.2$ and
$\dH(S^{\mathrm{grid}}_{3,\tau_3},S)=0.006$; the estimates are shown in
\Cref{fig:segment_deg_2}.

\section{The semicircle}
\label{app:semicircle}

\begin{figure}[h]
    \centering
    \includegraphics[width=0.7\linewidth]{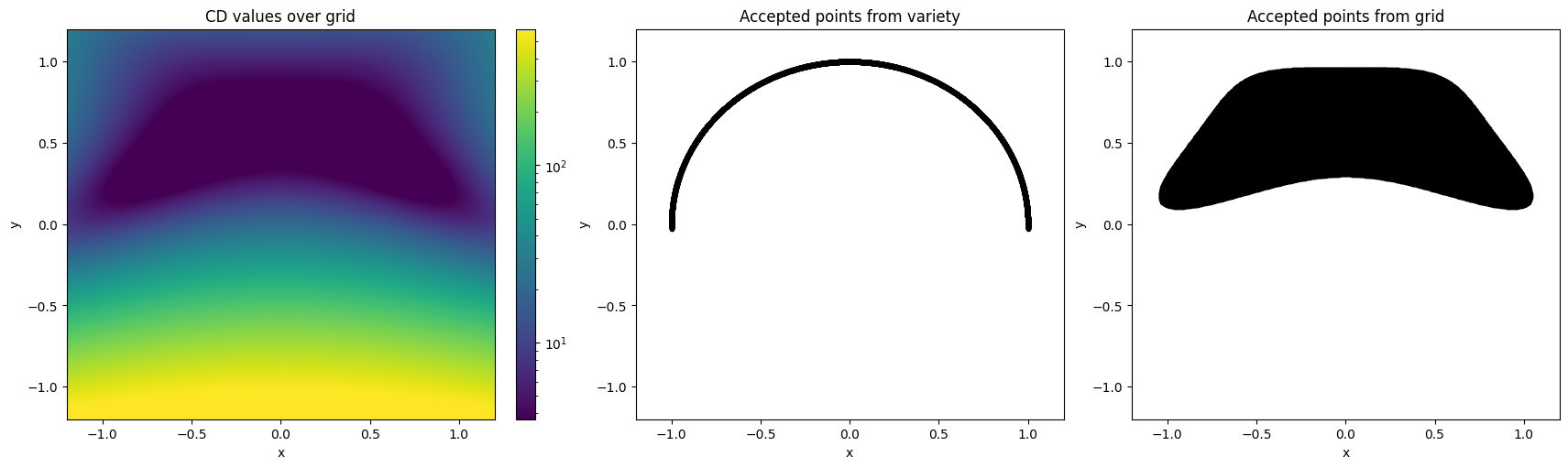}
    \caption{The CD kernel $\kappa_2$ for the semicircle. \emph{Left:} values
    of $\kappa_2$ on the $1000\times1000$ grid, on a logarithmic color scale;
    the kernel separates the two halves of the circle, although both have the
    same vanishing ideal. \emph{Middle:} the $1.1\times10^4$ accepted points of
    the variety, $\dH(S^{\mathrm{variety}}_{2,\tau_2},S)=0.05$. \emph{Right:}
    the $2\times10^5$ accepted grid points,
    $\dH(S^{\mathrm{grid}}_{2,\tau_2},S)=0.69$.}
    \label{fig:circle}
\end{figure}

\begin{figure}[h]
    \centering
    \includegraphics[width=0.7\linewidth]{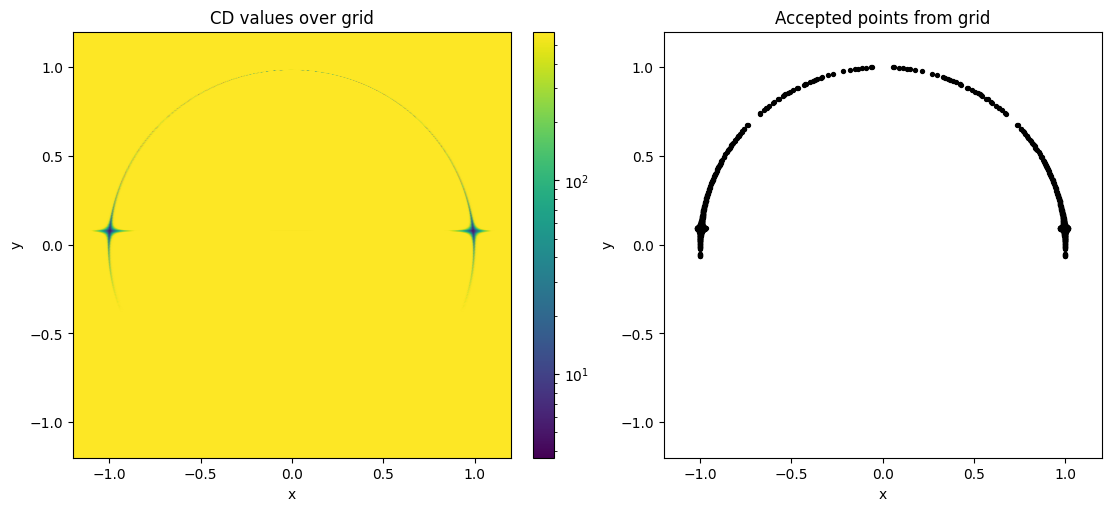}
    \caption{The CD kernel $\kappa_4$ for the semicircle. \emph{Left:} values
    of $\kappa_4$ on the $1000\times1000$ grid, on a logarithmic color scale.
    \emph{Right:} the accepted grid points,
    $\dH(S^{\mathrm{grid}}_{4,\tau_4},S)=0.19$.}
    \label{fig:circle_2}
\end{figure}

The polynomial optimization problem for the semicircle of
\Cref{sec:POP_setting}, which we solve at relaxation order~$2$, is
\begin{align*}
    \min_{x_1,x_2\in\R} \quad & F(x_1,x_2)=(x_1^2+x_2^2-1)^2 \\
    \text{s.t.}\quad & x_2\ge0,\quad 1-x_2\ge0 ,
\end{align*}
and its set of minimizers is
$S=\{(x_1,x_2)\in\R^2: x_1^2+x_2^2=1,\ x_2\ge0\}$. In degree $2$ the kernel of
the moment matrix contains the vector
$(-0.57734851,\,0,\,0,\,0.57734864,\,0,\,0.57735366)$, that is, the polynomial
$x_1^2+x_2^2-1$. The kernel can give no more than this, since the vanishing
ideal of the semicircle equals that of the \emph{whole} circle, and no
moment-matrix kernel in any degree will produce the constraint $x_2\ge0$. We
evaluate $\kappa_2$ at $2.2\times10^4$ samples on the circle and threshold at
$\tau_2^{\mathrm{variety}}=9.77$. The extreme points of the $1.1\times10^4$
accepted samples are $(-1.00000363,\,-1.7\times10^{-4})$ and
$(1.00000364,\,1.2\times10^{-4})$, and
$\dH(S^{\mathrm{variety}}_{2,\tau_2},S)=0.05$.
The lower half of the circle is rejected by the CD kernel alone, $\kappa_2$
increasing by more than an order of magnitude across $x_2=0$
(\Cref{fig:circle}, left). Thresholding the $10^6$ grid points in degree $2$
with $\tau_2^{\mathrm{grid}}=4.64$ gives
$\dH(S^{\mathrm{grid}}_{2,\tau_2},S)=0.69$, and thresholding in degree $4$ with
$\tau_4^{\mathrm{grid}}=30$ gives $\dH(S^{\mathrm{grid}}_{4,\tau_4},S)=0.19$
(\Cref{fig:circle_2}).

\section{The product of a semicircle and a segment}
\label{app:mixed}

\begin{figure}[h]
    \centering
    \includegraphics[width=0.8\linewidth]{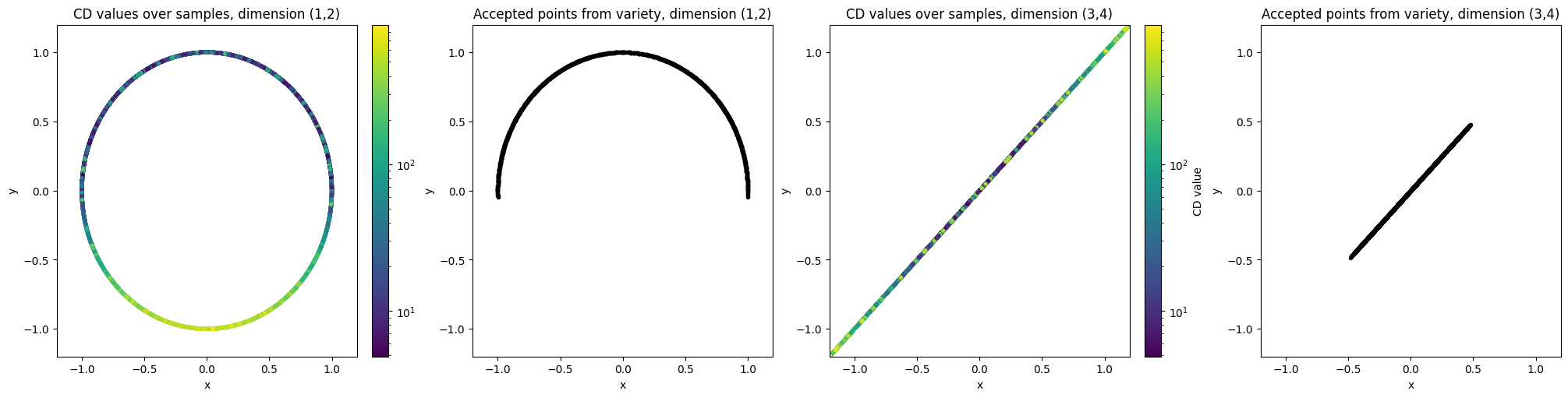}
    \caption{CD kernel $\kappa_2$ and the accepted variety points. \emph{Left:} $\kappa_2$ evaluated on the variety, projected onto dimensions 1 and 2. \emph{Middle left:} The accepted points projected onto dimensions 1 and 2. \emph{Middle right:} $\kappa_2$ evaluated on the variety, projected onto dimensions 3 and 4. \emph{Right:} The accepted points projected onto dimensions 3 and 4.}
    \label{fig:4d_variety}
\end{figure}

This example is set up in $\R^4$.
To include terms of higher degree in the objective while keeping the moment
matrix small, we take four variables and solve, at relaxation degree $4$,
\begin{align*}
    \min_{x\in\R^{4}} \quad & F(x)=(1-x_1^2-x_2^2)^2+(x_3-x_4)^2 \\
    \text{s.t.}\quad & x_2\ge0,\quad x_3\ge-0.5,\quad x_4\le0.5,
\end{align*}
whose set of minimizers is
\begin{equation*}
    S=\{(x_1,x_2)\in\R^2: x_1^2+x_2^2=1,\ x_2\ge0\}
      \times
      \{(x_3,x_4)\in\R^2: x_3=x_4,\ -0.5\le x_3\le0.5\}\subset\R^4 .
\end{equation*}
In degree $2$ the kernel of $M_2(y)$ is two-dimensional, spanned by the
coefficient vectors of $x_1^2+x_2^2-1$ and $x_3-x_4$. The homotopy-continuation
sampler, applied to these two polynomials, gives $10^4$ samples, at which we
evaluate the CD kernel $\kappa_2$ of degree $2$. We also evaluate the CD kernel $\kappa_4$ of degree
$4$ at $6\times10^6$ grid points.
Thresholding at $\tau_2^{\mathrm{variety}}= 14$ and $\tau_4^{\mathrm{grid}}= 20$ results into $S_{2,\tau_2}^{\mathrm{variety}}$ consisting of 2676 points and $S_{4,\tau_4}^{\mathrm{grid}}$ which has 324 points.
 We report the Hausdorff distances computed against a $0.01$-dense discretization of $S$:
\begin{equation*}
    \dH\big(S_{2,\tau_2}^{\mathrm{variety}},S\big)=0.1,
    \quad \dH\big(S_{4,\tau_4}^{\mathrm{grid}},S\big) = 0.246.
\end{equation*}
The projection of samples accepted on the variety is depicted in Figure \ref{fig:4d_variety}.

\section{Experimental details}
\label{app:details}

\paragraph{Examples with moments in closed form.}
For each of the five experiments of this kind in
\Cref{tab:recovery-experiments}, $\mu$ is the uniform measure with respect to
$\mathcal H^s|_{\mathcal Z}$ on the indicated subset of the indicated variety,
except for the second hyperbola, where $\mu$ is the image of the uniform
measure on $t\in[-T,T]$, $T=\operatorname{arcosh}2$, under
$t\mapsto(\cosh t,\sinh t)$; its density with respect to arc length is bounded
above and below, as \Cref{problem_main} requires. The moments up to degree $2d_v$ are computed in closed form. The carrier is
recovered from $\ker M_{d}(y)$ at the smallest $d$ for which the kernel is
nontrivial, which is $d=2$ for the conics, the cone and the cylinder.

\paragraph{Examples from polynomial optimization.}
The POPs for the segment and the semicircle in \Cref{sec:pop} are solved with
TSSOS \citep{magron2021tssos,wang2021tssos} at relaxation degree $2$, and the
resulting pseudo-moment sequence is used as $y$ in \Cref{alg:support} without
further processing. Neither pseudo-moment matrix has the flat extension
property, so that neither the extraction procedure of \citet{CurtoFialkow1998}
nor the GNS construction of \citet{LopezQuijorna2021DetectingGNS} is available.

\end{document}